\documentclass{article}
\usepackage{stmaryrd}
\usepackage{mathrsfs}
\usepackage[centertags]{amsmath}
\usepackage{amsfonts, dsfont}
\usepackage{amssymb}
\usepackage{amsthm}
\usepackage[OT2,T1]{fontenc}
\DeclareSymbolFont{cyrletters}{OT2}{wncyr}{m}{n}
\DeclareMathSymbol{\Sha}{\mathalpha}{cyrletters}{"58}

\input xy
\xyoption{all}

\theoremstyle{plain}
\newtheorem{thm}{Theorem}[section]
\newtheorem{cor}[thm]{Corollary}
\newtheorem{lem}[thm]{Lemma}
\newtheorem{prop}[thm]{Proposition}
\newtheorem*{thm2}{Theorem}

\theoremstyle{definition}
\newtheorem{defn}[thm]{Definition}
\newtheorem*{remark}{Remark}
\newtheorem*{ack}{Acknowledgments}

\newcommand{\bd}{\begin{defn}}
\newcommand{\ed}{\end{defn}}
\newcommand{\bl}{\begin{lem}}
\newcommand{\el}{\end{lem}}
\newcommand{\bp}{\begin{prop}}
\newcommand{\ep}{\end{prop}}
\newcommand{\bt}{\begin{thm}}
\newcommand{\et}{\end{thm}}
\newcommand{\bc}{\begin{cor}}
\newcommand{\ec}{\end{cor}}
\newcommand{\br}{\begin{remark}}
\newcommand{\er}{\end{remark}}
\newcommand{\bdi}{\begin{diagram}}
\newcommand{\edi}{\end{diagram}}
\newcommand{\beq}{\begin{equation}}
\newcommand{\eeq}{\end{equation}}
\newcommand{\ba}{\begin{array}}
\newcommand{\ea}{\end{array}}
\newcommand{\bpf}{\begin{proof}}
\newcommand{\epf}{\end{proof}}

\newcommand{\Z}{\mathds{Z}}
\newcommand{\Q}{\mathds{Q}}
\newcommand{\Zp}{\mathds{Z}_{p}}
\newcommand{\Qp}{\mathds{Q}_{p}}
\newcommand{\al}{\alpha}
\newcommand{\be}{\beta}
\newcommand{\Ga}{\Gamma}
\newcommand{\ga}{\gamma}

\newcommand{\la}{\lambda}

\DeclareMathOperator{\Sel}{Sel} \DeclareMathOperator{\Gal}{Gal}
 \DeclareMathOperator{\rank}{rank}
\DeclareMathOperator{\corank}{corank}

\newcommand{\cyc}{\mathrm{cyc}}
\newcommand{\m}{\mathfrak{m}}
\newcommand{\M}{\mathfrak{M}}
\newcommand{\p}{\mathfrak{p}}

\newcommand{\ot}{\otimes}
\newcommand{\ilim}{\displaystyle \mathop{\varinjlim}\limits}

\newcommand{\coker}{\mathrm{coker}\,}

\newcommand{\lra}{\longrightarrow}

\newcommand{\ps}[1]{\llbracket #1 \rrbracket}

\begin{document}
\title{On the growth of Mordell-Weil ranks in $p$-adic Lie extensions}
 \author{
 Pin-Chi Hung\footnote{Room R0718, First Academic Building, Soochow University, Taipei, R.O.C.
 E-mail: \texttt{pinchihung1111@gmail.com}} \quad
  Meng Fai Lim\footnote{School of Mathematics and Statistics $\&$ Hubei Key Laboratory of Mathematical Sciences,
Central China Normal University, Wuhan, 430079, P.R.China.
 E-mail: \texttt{limmf@mail.ccnu.edu.cn}} }
\date{}
\maketitle

\begin{abstract} \footnotesize
\noindent Let $p$ be an odd prime and $F_{\infty}$ a $p$-adic Lie extension of a number field $F$. Let $A$ be an abelian variety over $F$ which has ordinary reduction at every primes above $p$. Under various assumptions, we establish asymptotic upper bounds for the growth of Mordell-Weil rank of the abelian variety of $A$ in the said $p$-adic Lie extension. Our upper bound can be expressed in terms of invariants coming from the cyclotomic level. Motivated by this formula, we make a conjecture on an asymptotic upper bound of the growth of Mordell-Weil ranks over a $p$-adic Lie extension which is in terms of the Mordell-Weil rank of the abelian variety over the cyclotomic $\Zp$-extension. Finally, it is then natural to ask whether there is such a conjectural upper bound when the abelian variety has non-ordinary reduction. For this, we can at least modestly formulate an analogous conjectural upper bound for the growth of Mordell-Weil ranks of an elliptic curve with good supersingular reduction at the prime $p$ over a $\Zp^2$-extension of an imaginary quadratic field.


\end{abstract}

\section{Introduction}

Let $A$ be an Abelian variety defined over a number field $F$. The well-known Mordell-Weil theorem asserts that the group $A(F)$ of $F$-rational points is a finitely generated abelian group. In particular, this group has a well-defined $\Z$-rank which is called the Mordell-Weil rank of $A$. In this paper, we are interested in the variation of the Mordell-Weil ranks of an abelian variety in a $p$-adic Lie extension, where $p$ is an odd prime at which the abelian variety has ordinary reduction at every prime of $F$ above $p$.  In studying the Mordell-Weil rank, the Selmer group plays an important role. In his fundamental work \cite{Maz}, Mazur developed the (ordinary) Iwasawa theory of Selmer groups, and applied it to obtain an upper bound for the growth of Mordell-Weil ranks in a cyclotomic $\Zp$-extension which we now describe.

Let $F^{\cyc}$ be the cyclotomic $\Zp$-extension of $F$. Denote by $F_n$ the intermediate subfield of $F^{\cyc}/F$ with index $|F_n:F|=p^n$. Write $X(A/F^{\cyc})$ for the Pontryagin dual of the Selmer group of $A$ over $F^{\cyc}$. This Selmer group carries a natural $\Zp\ps{\Ga}$-module structure, where $\Ga=\Gal(F^{\cyc}/F)$. Mazur conjectured that $X(A/F^{\cyc})$ is a torsion $\Zp\ps{\Ga}$-module. Granted the validity of the conjecture, one can attach Iwasawa $\lambda$-invariant to this module which is usually denoted by $\la_{\Zp\ps{\Ga}}(X(A/F^{\cyc}))$. The following is a theorem of Mazur \cite[p. 185]{Maz} (or see \cite[Theorem 1.9]{GR99}) which gives a uniform bound on the Mordell-Weil ranks in a cyclotomic $\Zp$-extension in term of this Iwasawa $\lambda$-invariant.

\begin{thm2}[Mazur]
Let $A$ be an abelian variety defined over a number field which has good ordinary reduction at every primes above $p$. Let $F^{cyc}$ be the cyclotomic $\Zp$-extension of $F$ with intermediate subfield $F_n$ of index $|F_n:F|=p^n$. Suppose that $X(A/F^{\cyc})$ is a torsion $\Zp\ps{\Ga}$-module. Then we have
\[ \rank_{\Z}(A(F_n)) \leq \la_{\Zp\ps{\Ga}}(X(A/F^{\cyc})). \]
\end{thm2}

The goal of this paper is to search for such an analogous upper bound over a $p$-adic Lie extension of higher dimension. Indeed, if $F_{\infty}$ is now a uniform $p$-adic Lie extension (see Section 3 for definition) of $F$ of dimension $d$ with Galois group $G$,  there is a natural extension of the notion of a torsion $\Zp\ps{G}$-module (see \cite{V02}). Under the assumption that $X(A/F_{\infty})$ is torsion over $\Zp\ps{G}$, one can show that $\rank_{\Z}(A(F_n)) = O(p^{(d-1)n})$ by appealing to the work of Harris \cite{Har00} (also see \cite[Corollary 19]{Bh}, \cite[Corollary 2.9]{HV} or \cite[Theorem 3.2]{LS}). However, Harris's result does not give a concrete upper bound as in the cyclotomic $\Zp$-extension. The main reason behind this is that we do not have a nice enough structure theory for modules over noncommutative Iwasawa algebras unlike the cyclotomic situation (see \cite{CFKSV, CSSalg}).

After much intensive study by Coates, Fukaya, Kato, Sujatha and Venjakob \cite{CFKSV, V05}, they were led to conjecture that the dual Selmer group $X(A/F_{\infty})$ satisfies a stronger torsion property which enables one to define a higher analogue of the Iwasawa $\lambda$-invariant. We now describe this aspect of their work. Denote by $X(A/F_{\infty})(p)$ the $\Zp\ps{G}$-submodule of $X(A/F_{\infty})$ consisting of elements annihilated by some power of $p$ and write $X_f(A/F_{\infty}) = X(A/F_{\infty})/X(A/F_{\infty})(p)$. Coates et al \cite{CFKSV} conjectured that $X_f(A/F_{\infty})$ is finitely generated over $\Zp\ps{H}$, where $H=\Gal(F_{\infty}/F^{\cyc})$.  Granted this conjecture, it then makes sense to speak of $\rank_{\Zp\ps{H}}\big(X_f(A/F_{\infty})\big)$. It has been long observed in literature that this quantity serves as a higher analog of the classical $\la$-invariant (for instances, see \cite{CH, Ho}). In view of this, it would seem natural to expect an upper bound of the Mordell-Weil ranks which has a description in term of this quantity, and this is precisely the main theorem of our paper.

\begin{thm2}[Theorem \ref{main theorem}]
Assume that $(i)$ $A$ is an abelian variety over a number field $F$ which has ordinary reduction at every primes above $p$,
$(ii)$ $F_{\infty}$ is a uniform admissible $p$-adic extension of $F$ of dimension $d\geq 2$ and $(iii)$ $X_f(A/F_{\infty})$ is finitely
generated over $\Zp\ps{H}$. Denoting by $F_n$ the fixed field of $\Gal(F_{\infty}/F)^{p^n}$, we have
\[ \rank_{Z}(A(F_n))\leq \rank_{\Zp\ps{H}}\big(X_f(A/F_{\infty})\big)p^{(d-1)n} + O(p^{(d-2)n}).\]
\end{thm2}

We mention in passing that the error term $O(p^{(d-2)n})$ arises due to the usage of an asymptotic formula of Harris \cite{Har00}.
By imposing an extra assumption, we can elucidate the error terms further, and this is the content of the next theorem.

\begin{thm2}[Theorem \ref{main theorem2}]
Retain all the assumptions of Theorem \ref{main theorem}. Assume further that $H_i(H_n,X(A/F_{\infty}))$ is finite for every $i\geq 1$ and $n\geq 0$. Then one has
\[ \rank_{Z}(A(F_n))\leq \rank_{\Zp\ps{H}}\big(X_f(A/F_{\infty})\big)p^{(d-1)n} + d\corank_{\Zp}(A(F_{\infty})(p)).\]
\end{thm2}

The point of the extra finiteness assumption in the preceding theorem is to allow us to \textit{avoid} the usage of Harris's formula which is the key in obtaining such a precise upper bound. The finiteness assumption are known to be valid in many situations (see Remark after Theorem \ref{main theorem2}).

We should mention that a proof of Theorem \ref{main theorem2} was established in \cite[Corollary 1.4]{DL} by an algebraic $K$-theoretical argument. (Although their result is stated for (solvable) admissible $p$-adic Lie extension of dimension $\leq 3$, one can check that their algebraic $K$-theoretical argument carries over to the general situation.) Our proof here is different from there in that we do not use any algebraic $K$-theory, instead giving a direct proof via control theorems and some rank calculations of Howson and Harris. Our reason of having this approach is twofold.

Firstly, the above approach can be adapted to yield a description of $\rank_{\Zp\ps{H}}\big(X_f(A/F_{\infty})\big)$ in terms of invariants coming from the cyclotomic level $F^{\cyc}$. Combining this description with the above theorems, one obtains an upper bound in terms of these cyclotomic invariants (see Corollary \ref{rank formula corollary}). The latter inspires us to make a conjecture on an asymptotic upper bound of the Mordell-Weil ranks in terms of the cyclotomic Mordell-Weil rank (see Conjecture 1). We like to mention that although Theorems \ref{main theorem} and \ref{main theorem2}, and Corollary \ref{rank formula corollary} are derived under the validity of $\M_H(G)$-conjecture, our Conjecture 1 \textit{does not} require the $\M_H(G)$-conjecture in its formulation (although we need an appropriate conjecture of Mazur for our Conjecture 1). We provide some (mild) theoretical evidence to our Conjecture 1 (see Section \ref{Additional evidence for Conjecture 2}). We also mention that in proving Conjecture 1 in these situations, we do not assume the $\M_H(G)$-conjecture.

The second reason of adopting a non algebraic $K$-theoretical proof stems from a recent work of Lei and Sprung \cite{LeiS}, where they obtained an upper bound for an elliptic curve with good supersingular reduction at the prime $p$ over a $\Zp^2$-extension of an imaginary quadratic field which is in the spirit of Harris.  As the supersingular situation is slightly more delicate, the approach we adopted is more suitable than the algebraic K-theoretical approach. Indeed, we are able to establish analogue of Theorem \ref{main theorem2} for this non-ordinary situation under an appropriate supersingular variant of the $\M_H(G)$-conjecture. We also mention that as in the case of Theorem \ref{main theorem2}, we do not use Harris's asymptotic formula which allows us to establish a precise upper bound (see Theorem \ref{supersingular theorem} for details). Following the ordinary situation, we also formulate a conjecture on an upper bound of the Mordell-Weil ranks in this non-ordinary setting (see Conjecture 2).

We now give a brief description of the layout of the paper. In Section \ref{algebra},
we recall certain algebraic notion which will be used subsequently in the
paper. We also prove several lemmas in preparation for the proof of the main results.  Section \ref{Selmer} is where we introduce the Selmer groups and prove our main results. In Section \ref{rank section}, we calculate the quantity $\rank_{\Zp\ps{H}}\big(X_f(A/F_{\infty})\big)$ in terms of various cyclotomic invariants. It is also here that we state our Conjecture 1 and present some evidence for it. This is further continued in Section \ref{Additional evidence for Conjecture 2}, where we describe how the combination of the works of Cornut-Vatsal, Howard, Nekov\'{a}\u{r} can be applied to establish the validity of our Conjecture 1 for a $\Zp^2$-extension of an imaginary quadratic field. Finally, in Section \ref{elliptic supersingular conjecture}, we establish results analogue to those in Sections \ref{Selmer} and \ref{rank section} for an elliptic curve with good supersingular reduction over the $\Zp^2$-extension of an imaginary field. Building on this, we formulate our conjecture (Conjecture 2) on the upper bound of the Mordell-Weil ranks in this modest non-ordinary setting.

\begin{ack}
The authors are very grateful of Antonio Lei for his comments and interest on the paper. We would especially like to thank him for the discussion pertaining to Section \ref{elliptic supersingular conjecture} and for making us aware of the paper \cite{LV}. The authors also like to thank John Coates for his interest and helpful comments on the paper.  We also thank Ming-Lun Hsieh for his encouragement on the authors' collaboration. We are grateful to the anonymous referee
for providing various helpful comments and feedback.
Some part of the research of this article was conducted when M. F. Lim was visiting the National Center for Theoretical Sciences and Academia Sinica of Taiwan, and he would like to acknowledge the hospitality and conducive working conditions provided by these institutes. Finally, P. -C. Hung's research is supported by the MOST grant 107-2115-M-031-001-MY2, and M. F. Lim's research is supported by the National Natural Science Foundation of China under Grant No. 11550110172 and Grant No. 11771164.
 \end{ack}

\section{Algebraic Preliminaries} \label{algebra}

In this section, we recall some algebraic preliminaries that will be
required in the later part of the paper. Let $G$ be a compact pro-$p$ $p$-adic
Lie group without $p$-torsion. It is well known that $\Zp\ps{G}$ is
an Auslander regular ring (cf. \cite[Theorems 3.26]{V02}).
Furthermore, the ring $\Zp\ps{G}$ has no zero divisors (cf.\
\cite{Neu}), and therefore, admits a skew field $Q(G)$ which is flat
over $\Zp\ps{G}$ (see \cite[Chapters 6 and 10]{GW} or \cite[Chapter
4, \S 9 and \S 10]{Lam}). Thanks to this, we can define the notion of $\Zp\ps{G}$-rank of a finitely generated $\Zp\ps{G}$-module $M$, which is given by
$$ \rank_{\Zp\ps{G}}(M)  = \dim_{Q(G)} (Q(G)\ot_{\Zp\ps{G}}M). $$
The module $M$ is then said to be a
\textit{torsion} $\Zp\ps{G}$-module if $\rank_{\Zp\ps{G}} (M) = 0$.

Now if $M$ is a finitely generated $\Zp\ps{G}$-module, then its homology groups $H_i(G,M)$ are finitely generated over $\Zp$ (see \cite[Proof of Theorem 1.1]{Ho} or \cite[Lemma 3.2.3]{LS}). Hence the quantity $\rank_{\Zp}\big(H_i(G,M)\big)$ is well-defined. In view of this observation, we can now state the following result of Howson (see \cite[Theorem 1.1]{Ho} or \cite[Lemma 4.3]{LimFine}).

\bp[Howson] \label{Howson}
Let $M$ be a finitely generated $\Zp\ps{G}$-module. Then we have
\[ \rank_{\Zp\ps{G}}(M) = \sum_{i=0}^d(-1)^i\rank_{\Zp}\big(H_i(G,M)\big),\]
where here $d$ denotes the dimension of the $p$-adic group $G$.
\ep

From now on, our group $G$ is always assumed to be a uniform pro-$p$ group in the sense of \cite[Section 4]{DSMS}. We write $G_n$ for the lower
$p$-series $P_{n+1}(G)$ which is defined recursively by $P_{1}(G) = G$ and
\[ P_{n+1}(G) = \overline{P_{n}(G)^{p}[P_{n}(G),G]}, ~\mbox{for}~ n\geq 1. \]
It follows from \cite[Thm.\
3.6]{DSMS} that $G^{p^n} =
P_{n+1}(G)$ and that we have an equality $|G:P_2(G)| = |P_n(G):
P_{n+1}(G)|$ for every $i\geq 1$ (cf. \cite[Definition 4.1]{DSMS}). It follows from these that  $|G:G_n| = p^{dn}$, where $d= \dim G$. We now record the following lemma whose proof is left to the readers as an exercise (or see \cite[Corollary 1.5]{Ho}).

\bl \label{lemma rank}
 Let $M$ be a finitely generated $\Zp\ps{G}$-module. Then $M$ is finitely generated over $\Zp\ps{G_n}$ with
 \[\rank_{\Zp\ps{G_n}}(M) =|G:G_n|\rank_{\Zp\ps{G}}(M)=p^{dn}\rank_{\Zp\ps{G}}(M).\]
\el

The next lemma will be useful in the subsequent of the paper.

\bl \label{Howson refined}
Let $G$ be a uniform pro-$p$ group of dimension $d$ and $M$ a finitely generated $\Zp\ps{G}$-module. Suppose that $H_i(G_n,M)$ is finite for every $i\geq 1$ and $n\geq 0$. (Here $G_0$ is to be understood as $G$.) Then for every $n\geq 0$, we have
\[ \rank_{\Zp}(M_{G_n}) = \rank_{\Zp\ps{G_n}}(M) = p^{dn}\rank_{\Zp\ps{G}}(M).\]
\el

\bpf
This follows from combining Proposition \ref{Howson} and Lemma \ref{lemma rank}.
\epf

To prepare for the next lemma, we need to introduce some more notation. For a $\Zp\ps{G}$-module $M$, denote by $M(p)$ the $\Zp\ps{G}$-submodule of $M$ consisting of all the elements of $M$ which are annihilated by some power of $p$.

\bl \label{rank compare lemma}
Let $G$ be a uniform pro-$p$ group of dimension $d$ and $M$ a finitely generated $\Zp\ps{G}$-module. Then for every $i\geq 0$, we have
\[ \rank_{\Zp}\big(H_i(G, M)\big) =\rank_{\Zp}\big(H_i(G, M_f)\big),\]
where $M_f=M/M(p)$.
\el

\bpf
 From the short exact sequence
 \[ 0\lra M(p) \lra M\lra M_f\lra 0,\]
 we have an exact sequence
 \[ H_i\big(G,M(p)\big)\lra H_i(G,M)\lra H_i(G,M_f)\lra H_{i-1}\big(G,M(p)\big),\]
 where $H_{-1}\big(G,M(p)\big)$ is to be understood to be zero.
 Since the ring $\Zp\ps{G}$ is Noetherian and $M$ is finitely generated over $\Zp\ps{G}$, so is $M(p)$. Therefore, there exists a sufficiently large $t$ such that $p^t$ annihilates $M(p)$, and hence all its $G$-homology groups.  As the $G$-homology groups are finitely generated over $\Zp$ (see discussion before Proposition \ref{Howson}), they must therefore be finite. The equality of the lemma is now a consequence of this observation and the above four terms exact sequence. \epf

\bl \label{bound for fg Zp}
Let $G$ be a uniform pro-$p$ group of dimension $d$. Suppose that $M$ is a $\Zp\ps{G}$-module which is finitely generated over $\Zp$. Then we have
\[ \rank_{\Zp}\big(H_1(G, M)\big) \leq d\rank_{\Zp}(M).\]
\el

 \bpf
By virtue of Lemma \ref{rank compare lemma}, we may assume that $M$ is free as a $\Zp$-module. Under this said assumption, we have a short exact sequence
  \[ 0\lra  M\stackrel{\cdot p}{\lra} M\lra M/p\lra 0\]
  which in turn induces an injection $H_1(G,M)/p \hookrightarrow H_1(G,M/p)$.
 It then follows that
 \[ \rank_{\Zp}\big(H_1(G, M)\big) \leq \dim_{\mathbb{F}_p}\big(H_1(G,M)/p\big) \leq \dim_{\mathbb{F}_p}\big(H_1(G,M/p)\big). \]
 Now since $G$ is pro-$p$, the only simple discrete $G$-module is isomorphic to $\Z/p$ with a trival $G$-action (cf. \cite[Corollary 1.6.13]{NSW}). Hence we may apply a d\'{e}vissage argument to obtain the inequality
\[ \dim_{\mathbb{F}_p}\big(H_1(G,M/p)\big) \leq \dim_{\mathbb{F}_p}\big(H_1(G,\Z/p)\big)\dim_{\mathbb{F}_p}(M/p).\]
As $G$ is a uniform group of dimension $d$, we have $\dim_{\mathbb{F}_p}\big(H_1(G, \Z/p)\big)=d$ by \cite[Theorem 4.35]{DSMS}. Finally, as $M$ is assumed to be torsionfree, we have $\rank_{\Zp}(M) = \dim_{\mathbb{F}_p}(M/p)$. The required estimate is now a consequence of these observations.
\epf

 Suppose that the uniform group $G$ contains a closed normal subgroup $H$ with the property that $\Ga:=G/H\cong \Zp$. Since $\Ga$ is clearly a uniform group, it follows from \cite[Proposition 4.31]{DSMS} that $H$ is also a uniform  group. Write $H_n$ (resp., $\Ga_n$) for the lower $p$-series $P_{n+1}(H)$ of $H$ (resp., for $P_{n+1}(\Ga)$ of $\Ga$).
The next lemma records the relations between the lower $p$-series of the groups $H$, $G$ and $\Ga$.

\bl \label{group}
For every $n\geq 1$, we have $H_n = H\cap G_n$ and $G_n/H_n\cong\Ga_n$.
\el

\bpf
Since $H$ and $G$ are uniform, we have $H_n=H^{p^n}$ and $G_n = G^{p^n}$ (cf. \cite[Theorem 3.6]{DSMS}). Clearly, one has the inclusion $H^{p^n}\subseteq H\cap G^{p^n}$. Conversely, let $h\in H\cap G^{p^n}$. Then there exists $g\in G$ such that $h= g^{p^n}$ which in turn implies that the coset $gH$ is a torsion element in $G/H$. But since $G/H\cong\Zp$ has no $p$-torsion, this forces $g\in H$. Hence we have $h\in H^{p^n}$ and this proves the first equality. For the second equality, one simply observes that \[G_n/H_n = G^{p^n}/H^{p^n}\cong G^{p^n}H/H= (G/H)^{p^n}\cong \Ga^{p^n} = \Ga_n.\]
\epf

We end the section with one final lemma.

\bl \label{rank compare Mf}
Let $G$ be a uniform pro-$p$ group which contains a closed normal subgroup $H$ with the property that $\Ga:=G/H\cong \Zp$. Let $M$ be a finitely generated $\Zp\ps{G}$-module which has the properties that $M_f:=M/M(p)$ is finitely generated over $\Zp\ps{H}$ and that $H_i(H,M)$ is finitely generated over $\Zp$ for all $i\geq 1$. Then for every $i\geq 1$, we have
\[ \rank_{\Zp}\big(H_i(H, M)\big) =\rank_{\Zp}\big(H_i(H, M_f)\big).\]
\el

\bpf
 Taking $H$-homology
of the following short exact sequence
\[ 0\lra M(p)\lra M\lra M_f\lra 0,\]
we obtain a long exact sequence
\[  H_{i}\big(H, M(p)\big)\lra H_{i}(H,
M) \stackrel{f_i}\lra H_i(H, M_f)\lra H_{i-1}\big(H, M(p)\big)\]
for $i\geq 1$.
As seen in the proof of Lemma \ref{rank compare lemma} there exists a sufficiently large $t$ such that $p^t$ annihilates $M(p)$. It then follows from this that $p^t$ annihilates $H_i\big(H,M(p)\big)$. This in turn implies that $p^t$ annihilates $\ker f_i$ and $\mathrm{coker} f_i$.
But since $H_i(H,M)$, and therefore, $\ker f_i$ is finitely generated over $\Zp$, it follows that $\ker f_i$ is finite. Also, as $M_f$ is finitely generated over $\Zp\ps{H}$, the group $H_i(H, M_f)$ is therefore finitely generated over $\Zp$ which in turn implies the same holds for $\mathrm{coker} f_i$. But we have already seen that $\mathrm{coker} f_i$ is annihilated by $p^t$, and so $\mathrm{coker} f_i$ is finite. In conclusion, the map $f_i$ has finite kernel and cokernel, and the equality of the lemma is now an immediate consequence of this.
\epf

\section{Selmer groups} \label{Selmer}

In this section, we recall the definition of the
Selmer group of an abelian variety. As before, $p$ will denote an odd
prime. Let $F$ be a number field and $A$ an abelian variety
over $F$. Let $v$ be a prime of $F$. For every finite
extension $L$ of $F$, we define
 \[ J_v(A/L) = \bigoplus_{w|v}H^1(L_w, A)(p),\]
where $w$ runs over the (finite) set of primes of $L$ above $v$. If
$\mathcal{L}$ is an infinite extension of $F$, we define
\[ J_v(A/\mathcal{L}) = \ilim_L J_v(A/L),\]
where the direct limit is taken over all finite extensions $L$ of
$F$ contained in $\mathcal{L}$. For any algebraic (possibly
infinite) extension $\mathcal{L}$ of $F$, the Selmer group of $A$
over $\mathcal{L}$ is defined to be
\[ \Sel(A/\mathcal{L}) = \ker\Big(H^1(\mathcal{L}, A(p))\lra \bigoplus_{v} J_v(A/\mathcal{L})
\Big), \] where $v$ runs through all the primes of $F$.

If $L$ is a finite extension of $F$, then the Selmer group and the Mordell-Weil group are related by the following short exact sequence
\[  0\lra A(L)\ot_{\Z}\Qp/\Zp \lra \Sel(A/L) \lra \Sha(A/L)(p)\lra 0,     \]
where $\Sha(A/L)$ is the Tate-Shafarevich group. It follows from this that
\[ \rank_{\Z}(A(L))\leq \corank_{\Zp}\big(\Sel(A/L)\big).\]
Hence the problem of obtaining an upper bound for the Mordell-Weil ranks is reduced to obtaining an upper bound for the coranks of the Selmer groups.
(Of course, it has been conjectured that $\Sha(A/L)$ is finite and hence the above inequality should be an equality under this conjecture. However, for our purposes, we do not need to assume this.)

A Galois extension $F_{\infty}$ of $F$ is said to be a uniform admissible $p$-adic Lie
extension of $F$ if (i) $\Gal(F_{\infty}/F)$ is a uniform pro-$p$ group, (ii) $F_{\infty}$ contains the cyclotomic $\Zp$-extension
$F^{\cyc}$ of $F$ and (iii) $F_{\infty}$ is unramified outside a
finite set of primes of $F$. We shall always write $G = \Gal(F_{\infty}/F)$, $H =
\Gal(F_{\infty}/F^{\cyc})$ and $\Ga =\Gal(F^{\cyc}/F)$. Denote by $X(A/F_{\infty})$ the
Pontryagin dual of $\Sel(A/F_{\infty})$.

To continue, we need to recall certain facts from \cite{CG}. For now, let $K$ be a finite extension of $\Qp$. Write $I_K$ for the inertia subgroup. Suppose for now that $A$ is an abelian variety of dimension $g$ defined over $K$. Write $V=T_p(A)\ot_{\Zp}\Qp$, where $T_p(A)$ denotes the Tate module of $A$. Following \cite{CG}, we let $W$ be the $\Gal(\bar{K}/K)$-invariant $\Qp$-subspace of $V$ of minimal dimension such that some subgroup of $I_K$ of finite index acts trivially on $V/W$.
Set $C$ to be the image of $W$ under the natural map $V \longrightarrow V/T_p(A)=A(p)$. As seen in the discussion in \cite[pp 150]{CG}, in the event that the abelian variety $A$ has semistable reduction over $K$, $C$ is precisely $\mathcal{F}(\bar{\m})(p)$, where $\mathcal{F}$ is the formal group over $\mathcal{O}_K$ attached to the Neron model for $A$ over $\mathcal{O}_K$. Then as abelian groups, we have $C\cong (\Qp/\Zp)^h$, where $g\leq h \leq 2g$.
In general, an abelian variety $A$ will have semistable reduction over some finite extension $K'$ of $K$. If $\mathcal{F}'$ is the associated formal group over $K'$, then the above discussion yields $C=\mathcal{F}'(\bar{\m})(p)$. Hence we still have $C\cong (\Qp/\Zp)^h$ with $g\leq h \leq 2g$.

Returning to the global situation, we let $A$ be an abelian variety defined over a number field $F$.
For each prime $v$ of $F$ above $p$, denote by $C_v$ the $\Gal(\bar{F}_v/F_v)$-submodule of $A(p)$ with $h_v=\corank_{\Zp}(C_v)$ defined as in the preceding paragraph. We can now state the following conjecture.

\medskip \noindent \textbf{Conjecture (Mazur, Schneider).} $\rank_{\Zp\ps{G}}\big(X(A/F_{\infty})\big)=\displaystyle\sum_{v|p}(h_v-g)$.

\medskip
The conjecture was first stated by Mazur in \cite{Maz} for an abelian variety with good ordinary reduction over a cyclotomic $\Zp$-extension. This conjecture was extended to general abelian varieties by Schneider in \cite{Sch}.
For a general $p$-adic Lie extension, the conjecture was raised in \cite{OcV03} (also see \cite{CSSGL2, HO, HV}).

We shall say that our abelian variety has \textit{ordinary reduction} at all primes above $p$ if $h_v=g$ for every $v|p$. For instance, an elliptic curve which has either good ordinary reduction (in the usual sense) or multiplicative reduction at every prime above $p$ is ordinary in the above sense. Throughout this paper (with the exception of Section \ref{elliptic supersingular conjecture}), we always assume that our abelian variety has \textit{ordinary reduction} at all primes above $p$. Under this assumption, the conjecture of Mazur and Schneider is then equivalent to saying that $X(A/F_{\infty})$ is a torsion $\Zp\ps{G}$-module.

Coates et al \cite{CFKSV, CS12} have predicted that $X(A/F_{\infty})$ satisfies a stronger torsion property. In fact, they formulated their conjecture on the structure of the dual Selmer group of an elliptic curve with good ordinary reduction at all primes above $p$. When the elliptic curve has multiplicative reduction at primes above $p$, this was formulated in \cite{Lee}. Here we merely mimic these prior works in stating the following conjecture for abelian variety with ordinary reduction (in the above sense) at all primes of $F$ above $p$.

\medskip \noindent \textbf{$\M_H(G)$-Conjecture.} \textit{For every admissible $p$-adic Lie extension $F_{\infty}$
of $F$,  $X(A/F_{\infty})/X(A/F_{\infty})(p)$ is a finitely
generated $\Zp\ps{H}$-module.}

\medskip
From now on, we write $X_f(A/F_{\infty}) =
X(A/F_{\infty})/X(A/F_{\infty})(p)$. Assuming the validity of the $\M_H(G)$-Conjecture, it then makes sense to speak of $\rank_{\Zp\ps{H}}\big(X_f(A/F_{\infty})\big)$. We can now state the following theorems.

\bt \label{main theorem}
Assume that $(i)$ $A$ is an abelian variety over a number field $F$ which has ordinary reduction at every prime above $p$,
$(ii)$ $F_{\infty}$ is a uniform admissible $p$-adic extension of $F$ of dimension $d\geq 2$ and $(iii)$ $X_f(A/F_{\infty})$ is finitely
generated over $\Zp\ps{H}$. Then we have
\[ \rank_{\Z}(A(F_n))\leq \rank_{\Zp\ps{H}}\big(X_f(A/F_{\infty})\big)p^{(d-1)n} + O(p^{(d-2)n}).\]
\et

As mentioned in the introduction, we can obtain a more precise upper bound under an additional assumption on the $H_n$-homology of the dual Selmer groups.

\bt \label{main theorem2}
Retain all the assumptions of Theorem \ref{main theorem}. Assume further that $H_i(H_n,X(A/F_{\infty}))$ is finite for every $i\geq 1$ and $n\geq 0$. Then
\[ \rank_{\Z}(A(F_n))\leq \rank_{\Zp\ps{H}}\big(X_f(A/F_{\infty})\big)p^{(d-1)n} + d\corank_{\Zp}\big(A(F_{\infty})(p)\big).\]
\et

\br
\begin{enumerate}
 \item[$(1)$] We have presented our results for uniform $p$-adic Lie extension mainly for convenience. By virtue of Lazard's theorem (see \cite[Corollary 8.34]{DSMS}), a compact $p$-adic Lie group contains a uniform subgroup of finite index. Therefore, by base-changing of the base field, we can obtain an upper bound of the Mordell-Weil ranks in an arbitrary admissible $p$-adic Lie extension.

      \item[$(2)$] As already mentioned in the introduction, Theorem \ref{main theorem2} can be proved by an algebraic $K$-theoretical argument similar to that in \cite[Corollary 1.4]{DL}. Our proof here is different from there in that we do not use any algebraic $K$-theory. A version of Theorem \ref{main theorem2} was also obtained in \cite[Proposition 6.9]{CH} for an elliptic curve over an $GL_2$-extension under the stronger assumption that $X(A/F_{\infty})$ is finitely generated over $\Zp\ps{H}$.

\item[$(3)$]  The extra assumption on the finiteness of the $H_n$-homology groups in the preceding theorem is known to be satisfied in many cases and we shall mention them here.
\begin{enumerate}
\item[$(a)$] When $\dim G=2$ (i.e., $\dim H=1$), it follows from \cite[Proposition 5.1(c)]{LimMHG} that $H_i(H_n,X(A/F_{\infty}))=0$ for every $i\geq 1$ and $n\geq 0$. Hence this assumption holds in this situation.

\item[$(b)$] In fact if $\dim G\leq 3$, this assumption is also verified in \cite[Lemma 2.3]{DL} and is an intermediate argument required for the proof of \cite[Corollary 1.4]{DL}.

\item[$(c)$] If $A$ is an elliptic curve, the hypothesis has been verified for a large class of $p$-adic Lie extensions (see \cite[Proposition 13]{Bh}, \cite[Remark 2.6]{CSSalg} and \cite[Theorem 1.2 and Lemma 4.3]{Ze}).
\end{enumerate}
 \end{enumerate}
\er

Before proving the theorems, we first establish the following lemma.

\bl \label{main lemma}
Assume that $(i)$ $A$ is an abelian variety over a number field $F$ which has ordinary reduction at every prime above $p$ and
$(ii)$ $F_{\infty}$ is a uniform admissible $p$-adic extension of $F$. Then we have
\[ \rank_{\Z}(A(F_n))\leq \rank_{\Zp}\big(X_f(A/F_{\infty})_{G_n}\big) + \corank_{\Zp}\Big(H^1\big(G_n,A(F_{\infty})(p)\big)\Big).\]
\el

\bpf
It suffices to show that the quantity on the right is an upper bound for $\corank_{\Zp}(\Sel(A/F_n))$.
 Consider the following commutative diagram
\[ \xymatrixrowsep{0.4in}
\xymatrixcolsep{0.2in} \entrymodifiers={!! <0pt, .8ex>+}
\SelectTips{eu}{}\xymatrix{
    0 \ar[r]^{} & \Sel(A/F_n) \ar[d]^{s_n} \ar[r] &  H^1(G_S(F_n), A(p)) \ar[d]^{h_n}
    \ar[r] & \bigoplus_{v\in S}J_v(A/F_n) \ar[d]^{g_n} \\
    0 \ar[r]^{} & \Sel(A/F_{\infty})^{G_n} \ar[r] & H^1(G_S(F_{\infty}), A(p))^{G_n} \ar[r] &
    \Big(\bigoplus_{v\in S}J_v(A/F_{\infty})\Big)^{G_n}  } \]
with exact rows, and where the vertical maps are given by restriction maps. A diagram chasing argument immediately yields an exact sequence
\[ 0\lra \ker s_n \lra S(A/F_n) \lra S(A/F_{\infty})^{G_n}\] with $\ker s_n$ contained in $H^1(G_n,A(F_{\infty})[p^{\infty}])$. It follows from this that
\[ \corank_{\Zp}(S(A/F_n))\leq \rank_{\Zp}\big(X(A/F_{\infty})_{G_n}\big) + \corank_{\Zp}\Big(H^1\big(G_n,A(F_{\infty})(p)\big)\Big).\]
Finally, taking into account that $\rank_{\Zp}\big(X(A/F_{\infty})_{G_n}\big) = \rank_{\Zp}\big(X_f(A/F_{\infty})_{G_n}\big)$ by Lemma \ref{rank compare lemma}, we have the lemma.
\epf

We are now in position to prove our theorems.

\bpf[Proof of Theorem \ref{main theorem}]
  It follows from Lemma \ref{main lemma} that \[ \rank_{\Z}(A(F_n))\leq \rank_{\Zp}\big(X_f(A/F_{\infty})_{G_n}\big) + \corank_{\Zp}\Big(H^1\big(G_n,A(F_{\infty})(p)\big)\Big).\]
  By Lemma \ref{bound for fg Zp}, the second quantity on the right is bounded by $d\corank_{\Zp}\big(A(F_{\infty})(p)\big)$. It therefore remains to estimate $\rank_{\Zp}\big(X_f(A/F_{\infty})_{G_n}\big)$. Since $X_f(A/F_{\infty})$ is finitely generated over $\Zp\ps{H}$, it is also finitely generated over $\Zp\ps{H_n}$. Hence we have
  \[ \ba{c} \rank_{\Zp}\big(X_f(A/F_{\infty})_{G_n}\big) = \rank_{\Zp}\big(\big(X_f(A/F_{\infty})_{H_n}\big)_{\Ga_n}\big)\leq \rank_{\Zp}\big(X_f(A/F_{\infty})_{H_n}\big)\\
   =  \rank_{\Zp\ps{H}}\big(X_f(A/F_{\infty})\big)p^{(d-1)n} + O(p^{(d-2)n}),
  \ea
  \]
  where the first equality follows from Lemma \ref{group}, and the final equality follows from \cite[Theorem 1.10]{Har00} and noting that $H$ has dimension $d-1$. The required estimate is now immediate from combining the above estimates.
\epf

\bpf[Proof of Theorem \ref{main theorem2}]
As seen in the proof of Theorem \ref{main theorem}, we have
\[ \rank_{\Z}(A(F_n))\leq \rank_{\Zp}\big(X_f(A/F_{\infty})_{H_n}\big) + \corank_{\Zp}\Big(H^1\big(G_n,A(F_{\infty})(p)\big)\Big).\]
By Lemma \ref{bound for fg Zp}, the second quantity is bounded by $d\corank_{\Zp}(A(F_{\infty})(p))$.
On the other hand, it follows from Lemma \ref{rank compare lemma} and the hypothesis of the theorem that $H_i(H_n,X_f(A/F_{\infty}))$ is finite for every $i\geq 1$ and $n\geq 0$. Hence we may apply Lemma \ref{Howson refined} to conclude that
\[ \rank_{\Zp}\big(X_f(A/F_{\infty})_{H_n}\big) = \rank_{\Zp\ps{H}}\big(X_f(A/F_{\infty})\big)p^{(d-1)n}.\]
The conclusion of the theorem then follows from these.
\epf

\section{Comparing ranks} \label{rank section}

Retain the setting and notation from the previous section. We shall derive a formula which relates the $\rank_{\Zp\ps{H}}(X_f(A/F_{\infty}))$ in terms of certain invariants from the cyclotomic level.
Recall that by \cite[p. 150-151]{CG} (also see discussion in Section 3), for each prime $v$ of $F$ above $p$, we have a short exact sequence
\[ 0\lra C_v\lra A(p) \lra D_v\lra 0 \]
of discrete $\Gal(\bar{F}_v/F_v)$-modules which is characterized by the facts that $C_v$ is divisible and $D_v$ is the maximal quotient of $A(p)$ by a divisible subgroup such that the inertia group acts on $D_v$ via a finite quotient.
 Since our abelian variety $A$ has ordinary reduction at the prime $v$, both $C_v$ and $D_v$ are divisible abelian groups of corank $\dim(A)$. Furthermore, by \cite[Propositions 4.1, 4.7 and 4.8]{CG}, we have
\[ J_v(A/F_{\infty}) \cong \begin{cases}   \ilim_\mathcal{L}  \bigoplus_{w|v}H^1(\mathcal{L}_w, D_v),& \mbox{if } v\mbox{ divides }p \\
      \ilim_\mathcal{L}  \bigoplus_{w|v}H^1(\mathcal{L}_w, A(p)), & \mbox{if } v\mbox{ does not divides }p \end{cases}
 \] where the direct limit is taken over all finite extensions
$\mathcal{L}$ of $F^{\cyc}$ contained in $F_{\infty}$.

We can now state the main result of this section.

\bp \label{rank formula}
Let $F_{\infty}$ be a strongly admissible pro-$p$ Lie extension of $F$. Let $A$ be an abelian variety over $F$ which has ordinary reduction at every primes above $p$. Assume that $X(A/F_{\infty})$ satisfies the $\M_H(G)$-conjecture. Then
\[ \rank_{\Zp\ps{H}}\big(X_f(A/F_{\infty})\big) = \rank_{\Zp}\big(X_f(A/F^{\cyc})\big) +
\sum_{\substack{w\in S(F^{\cyc}),
\\ \dim H_w\geq 1}}\corank_{\Zp}\big(Z_v(F^{\cyc}_w)(p)\big),\]
where $S(F^{\cyc})$ is the set of primes of $F^{\cyc}$ above $S$. Here $Z_v$ denotes $D_v$ or $A(p)$ accordingly as $v$ divides $p$ or not.
\ep

\br
Proposition \ref{rank formula} has been proved for an elliptic curve $E$ under the stronger assumption that $X(E/F_{\infty})$ is finitely generated over $\Zp\ps{H}$ (see \cite[Theorem 16]{Bh},
\cite[Corollary 6.10]{CH}, \cite[Theorem 5.4]{HS}, \cite[Theorem
3.1]{HV} and \cite[Theorem 2.8]{Ho}). The approach for the proof in this general case follows those in the above citations. For the convenience of the readers, we shall supply a proof here.
\er

As a start, we have the following lemma.

\bl \label{short exact sequences}
Retaining the assumptions of Proposition \ref{rank formula}, we have short exact sequences
\[ 0 \lra \Sel(A/F^{\cyc})\lra H^1(G_S(F^{\cyc}), A(p))  \lra \bigoplus_{v\in S} J_v(A/F^{\cyc})\lra 0\]
and
\[ 0 \lra \Sel(A/F_{\infty})\lra H^1(G_S(F_{\infty}), A(p))  \lra \bigoplus_{v\in S} J_v(A/F_{\infty})\lra 0.\]
\el

\bpf
Since $X(A/F_{\infty})$ is assumed to satisfy the $\M_H(G)$-conjecture, it follows from \cite[Proposition 2.5]{CS12} that for every finite extension $L$ of $F$ contained in $F_{\infty}$, $X(A/L^{\cyc})$ is torsion over $\Zp\ps{\Ga_L}$, where $\Ga_L=\Gal(L^{\cyc}/L)$. Since $A(L^{\cyc})(p)$ is finite (cf. \cite{Ri}), we may apply a similar argument to that in \cite[Proposition 3.3]{LimMHG} to obtain a short exact sequence
\[ 0 \lra \Sel(A/L^{\cyc})\lra H^1(G_S(L^{\cyc}), A(p))     \lra \bigoplus_{v\in S} J_v(A/L^{\cyc})\lra 0.\]
In particular, this yields the first short exact sequence by taking $L=F$.
On the other hand, by taking direct limit over $L$, we obtain the second short exact sequence.
\epf

The next two lemmas are concerned with the $H$-homology of global cohomology groups and local cohomology groups.

\bl \label{global calculation}
Retain the assumptions of Proposition \ref{rank formula}. We then have that $H^i(H, H^1(G_S(F_{\infty}),A(p)))$ is cofinitely generated over
$\Zp$ for every $i\geq 1$. Moreover, we have an exact sequence
\[ 0 \lra H^1(H,
A(F_{\infty})(p))\lra H^1(G_S(F^{\cyc}), A(p))     \lra H^1(G_S(F_{\infty}), A(p))^H\lra H^2(H,
A(F_{\infty})(p))\lra 0\]
and isomorphisms \[
H^i\big(H, H^1(G_S(F_{\infty}), A(p))\big)\cong H^{i+2}\big(H,
A(F_{\infty})(p)\big) \mbox{ for } i\geq 1.\]
\el

\bpf
Since $X(A/F_{\infty})$ is assumed to satisfy the $\M_H(G)$-conjecture, it follows from \cite[Proposition 2.5]{CS12} that for every finite extension $L$ of $F$ contained in $F_{\infty}$, $X(A/L^{\cyc})$ is torsion over $\Zp\ps{\Ga_L}$, where $\Ga_L=\Gal(L^{\cyc}/L)$. Via similar arguments to those in \cite[Proposition 3.3 and Corollary 3.4]{LimMHG}, we have that
$H^2(G_S(F^{\cyc}),A(p))=0$ and $H^2(G_S(F_{\infty}),A(p))=0$. Hence the spectral sequence
\[ H^i\big(H, H^j(G_S(F_{\infty}), A(p))\big)\Longrightarrow
H^{i+j}(G_S(F^{\cyc}), A(p))\] degenerates to yield an exact sequence
\[ 0 \lra H^1(H,
A(F_{\infty})(p))\lra H^1(G_S(F^{\cyc}), A(p))     \lra H^1(G_S(F_{\infty}), A(p))^H\lra H^2(H,
A(F_{\infty})(p))\lra 0\]
and isomorphisms \[
H^i\big(H, H^1(G_S(F_{\infty}), A(p))\big)\cong H^{i+2}\big(H,
A(F_{\infty})(p)\big) \mbox{ for } i\geq 1\]
where the $\Zp$-cofinitely generation of latter
groups follow on noting that for any $p$-adic Lie group $H$ and any $\Zp$-cofinitely generated
$H$-module $W$, all of the cohomology groups $H^i(H,W)$ are cofinitely generated
$\Zp$-modules.
\epf

\bl \label{local calculation}
Retain the assumption of Proposition \ref{rank formula}. Then  $H^i\left(H, \bigoplus_{v\in S}J_v(A/F_{\infty})\right)$ is cofinitely generated over
$\Zp$ for every $i\geq 1$. Moreover, we have an exact sequence
\[ 0 \lra \bigoplus_{w\in S(F^{\cyc})}H^1(H_w,
Z_v(F_{\infty,w}))\lra \bigoplus_{v\in S}J_v(A/F^{\cyc})   \lra \left(
\bigoplus_{v\in S}J_v(A/F_{\infty})\right)^H \] \[\lra  \bigoplus_{w\in S(F^{\cyc})}H^2(H_w,
Z_v(F_{\infty,w}))\lra 0\]
and isomorphisms
\[
H^i\left(H, \bigoplus_{v\in S}J_v(A/F_{\infty})\right)\cong\bigoplus_{w\in S(F^{\cyc})}H^{i+2}\big(H_w,
Z_v(F_{\infty,w})\big) \mbox{ for } i\geq 1.\]
Here $Z_v$ denotes $D_v$ or $A(p)$ accordingly as $v$ divides $p$ or not.
\el

\bpf
This is a local version of Lemma \ref{global calculation} with a similar proof noting that
$H^2(F^{\cyc}_w,A(p))=0$ and $H^2(F_{\infty,w},A(p))=0$ by \cite[Theorem 7.1.8(i)]{NSW}.
\epf

We can now give the proof of Proposition \ref{rank formula}.

\bpf[Proof of Proposition \ref{rank formula}]
 Consider the following commutative diagram
\[  \entrymodifiers={!! <0pt, .8ex>+} \SelectTips{eu}{}\xymatrix{
    0 \ar[r]^{} & \Sel(A/F^{\cyc}) \ar[d]_{\al} \ar[r] &  H^1(G_S(F^{\cyc}), A(p))
    \ar[d]_{\be}
    \ar[r] & \displaystyle\bigoplus_{v\in S}J_v(A/F^{\cyc}) \ar[d]_{\ga} \ar[r] & 0 &\\
    0 \ar[r]^{} & \Sel(A/F_{\infty})^H \ar[r]^{} & H^1(G_S(F_{\infty}), A(p))^H \ar[r] & \
    \displaystyle\bigoplus_{v\in S}J_v(A/F_{\infty})^H \ar[r] &  H^1\big(H, S(A/F_{\infty})\big) \ar[r] & \cdots } \]
with exact rows.  To simplify notation, we write $W_{\infty}=
H^1(G_S(F_{\infty}), A(p))$ and
 $J_{\infty} = \displaystyle\bigoplus_{v\in S}J_v(A/F_{\infty})$. From the commutative diagram, we have a long exact
 sequence
 \[ \ba{c} 0\lra \ker\al \lra \ker \be
 \lra \ker \ga
 \lra \coker \al \lra \coker \be \\
  \lra \coker\ga \lra H^1\big(H, \Sel(A/F_{\infty})\big)
 \lra H^1(H, W_{\infty})
 \lra H^1(H, J_{\infty})\lra \cdots \\
 \cdots\lra H^{i-1}(H, J_{\infty}) \lra H^i\big(H, \Sel(A/F_{\infty})\big)
 \lra H^i(H, W_{\infty})
 \lra H^i(H, J_{\infty})\lra \cdots .\ea \]

By Lemmas \ref{global calculation} and \ref{local calculation}, the groups $\ker \be$, $\ker \ga$, $\coker \be$, $\coker\ga$, $ H^i(H, W_{\infty})$ (for $i\geq 1$) and $H^{i}(H, J_{\infty})$ (for $i\geq 1$) are cofinitely generated over $\Zp$. Thus, combining this observation with the above exact sequence, we have that $\ker \al$, $\coker \al$ and $H^1\big(H, S(A/F_{\infty})\big)$ (for $i\geq 1)$ are cofinitely generated over $\Zp$. Moreover, we have
\[ \ba{c} \corank_{\Zp}(\ker\al) - \corank_{\Zp}(\coker\al)=  -\displaystyle\sum_{i\geq 1}(-1)^i\corank_{\Zp}H^i(H,\Sel(A/F_{\infty}))\hspace{2in} \\
   \hspace{1in}+\displaystyle\sum_{i\geq 1}(-1)^{i}\corank_{\Zp}H^i(H,
A(F_{\infty})(p))-\displaystyle\sum_{\substack{w\in S(F^{\cyc}),
\\ \dim H_w\geq 1}}\left(\sum_{i\geq 1}(-1)^{i}\corank_{\Zp}H^i(H_w,
Z_v(F_{\infty,w}))\right), \ea \]
where here $Z_v$ denotes $D_v$ or $A(p)$ accordingly as $v$ divides $p$ or not.
Applying Proposition \ref{Howson} and Lemma \ref{rank compare Mf}, the right hand side is just
\[ -\displaystyle\sum_{i\geq 1}(-1)^i\rank_{\Zp}H_i(H,X_f(A/F_{\infty}))-\corank_{\Zp}H^0(H,
A(F_{\infty})(p))
 +\displaystyle\sum_{\substack{w\in S(F^{\cyc}),
\\ \dim H_w\geq 1}}\corank_{\Zp}H^0(H_w,
Z_v(F_{\infty,w})).
\]
Now consider the following commutative diagram
\[  \entrymodifiers={!! <0pt, .8ex>+} \SelectTips{eu}{}\xymatrix{
     & X(A/F_{\infty})(p)_H \ar[d]_{h'} \ar[r] &  X(A/F_{\infty})_H
    \ar[d]_{\al^{\vee}}
    \ar[r] & X_f(A/F_{\infty})_H \ar[d]_{h''} \ar[r]& 0 \\
    0 \ar[r]^{} & X(A/F^{\cyc})(p) \ar[r]^{} &  X(A/F^{\cyc}) \ar[r] & \
     X_f(A/F^{\cyc}) \ar[r] &  0 } \]
with exact rows. This in turns yields a long exact
 sequence
 \[   \ker h' \lra \ker (\al^{\vee})
 \stackrel{f}{\lra} \ker h''
 \lra \coker h' \lra \coker (\al^{\vee}) \lra \ker h''\lra 0.  \]
Since $X(A/F_{\infty})$ satisfies $\M_H(G)$-conjecture, we have that $X_f(A/F_{\infty})_H$ is finitely generated over $\Zp$. By \cite[Proposition 2.5]{CS12}, $X(A/F^{\cyc})$ is torsion over $\Zp\ps{\Ga}$ and so $X_f(A/F^{\cyc})$ is finitely generated over $\Zp$. Hence $\ker  h''$ and $\coker h''$ are finitely generated over $\Zp$, and we have
\[ \rank_{\Zp} (\ker h'') - \rank_{\Zp}(\coker h'')= \rank_{\Zp}\big(X_f(A/F_{\infty})_H\big) -\rank_{\Zp}(X_f(A/F^{\cyc}))\]

On the other hand, as already seen above, $\ker (\al^{\vee})$ and $\coker (\al^{\vee})$ are finitely generated over $\Zp$. Hence so are $\ker f$ and $\coker h'$. But since these latter groups are $p$-primary, they must be finite. Thus, we have
\[ \rank_{\Zp} (\ker (\al^{\vee})) - \rank_{\Zp}(\coker (\al^{\vee}))= \rank_{\Zp}\big(X_f(A/F_{\infty})_H\big) -\rank_{\Zp}(X_f(A/F^{\cyc})).\]
 Combining this with the above calculations and applying Proposition \ref{Howson} for $X_f(A/F_{\infty})$, we obtain the required formula.
\epf

 A combination of Theorem \ref{main theorem}/\ref{main theorem2} and Proposition \ref{rank formula} yields the following.

\bc \label{rank formula corollary}
Retain the setting of Theorem \ref{main theorem}.  Then we have
\[ \rank_{\Z}(A(F_n)) \leq \left(\rank_{\Zp}(X_f(A/F^{\cyc})) +
\sum_{\substack{w\in S(F^{\cyc}),
\\ \dim H_w\geq 1}}\corank_{\Zp}\big(Z_v(F^{\cyc}_w)(p)\big)\right)p^{(d-1)n} + O(p^{(d-2)n}),\]
where $S(F^{\cyc})$ is the set of primes of $F^{\cyc}$ above $S$.

Furthermore, in the event that the extra assumption of Theorem \ref{main theorem2} is also valid, we can replace $O(p^{(d-2)n})$ in the above inequality by $d\corank_{\Zp}\big(A(F_{\infty})(p)\big)$.
\ec

Notice that the bound on the right hand side makes sense as long as we know that $X(A/F^{\cyc})$ is torsion over $\Zp\ps{\Ga}$. Hence we are naturally led to raise the following question.

\medskip\noindent\textbf{Question.}
Can one prove the inequality in \ref{rank formula corollary} under the weaker assumption that $X(A/F^{\cyc})$ is torsion over $\Zp\ps{\Ga}$?
\medskip

Note that it is still unknown if one can deduce the validity of $\M_H(G)$-conjecture from Mazur's conjecture. A case where such an implication holds is when $X(A/F^{\cyc})$ is finitely
generated over $\Zp$ (see \cite[Theorem 2.1]{CS12}). However, there are examples where $X(A/F^{\cyc})$ is not finitely generated over $\Zp$ (see \cite{Dr} and \cite[\S1, Example 2]{Maz}). Despite this, it seems reasonable to conjecture that the upper bound is valid whenever $X(A/F^{\cyc})$ is torsion over $\Zp\ps{\Ga}$. In fact, we shall go one step further in formulating a more refined conjectural upper bound (see Conjecture 1 below). Before doing so, we recall the following observation of Mazur \cite{Maz}.

\bl \label{Mazur observation}
 Let $A$ be an abelian variety over a number field $F$ which has ordinary reduction at every primes above $p$. Suppose that $X(A/F^{\cyc})$ is finitely
generated torsion over $\Zp\ps{\Ga}$. Then $A(F^{\cyc})$ is finitely generated as an abelian group.
\el

\bpf
 Since Mazur's conjecture holds, we may apply his theorem (as mentioned in the introduction) to see that $\rank_{\Z}(A(F_n)) \leq \la_{\Zp\ps{\Ga}}(X(A/F^{\cyc}))$ for every $n$. Choose $n_0$ such that $\rank_{\Z}(A(F_{n_0}))$ is as large as possible. Then $A(F^{\cyc})/A(F_{n_0})$ must be a torsion group. Let $P\in A(F^{\cyc})$. Then there exists an integer $m\geq 1$ such that $mP\in A(F_{n_0})$. This in turn implies that $m\left(\sigma(P)-P\right)= \sigma(mP)-mP = 0$ for all $\sigma\in\Gal(F^{\cyc}/F_{n_0})$. In other words, $\sigma(P)-P\in A(F^{\cyc})_{tor}$. By a result of Ribet \cite{Ri}, the torsion subgroup $A(F^{\cyc})_{tor}$ of $A(F^{\cyc})$ is finite. Set $t= |A(F^{\cyc})_{tor}|$. Then $t\left(\sigma(P)-P\right) = 0$ or $\sigma(tP) = tP$ for all $\sigma\in\Gal(F^{\cyc}/F_{n_0})$. Hence $tP\in A(F_{n_0})$.

 Therefore, we can define a homomorphism $\varphi: A(F^{\cyc})\lra A(F_{n_0})$ by $P\mapsto tP$.
 The image of $\varphi$ is finitely generated since it is a subgroup of $A(F_{n_0})$. On the other hand, the kernel of $\varphi$ is $A(F^{\cyc})_{tor}$ and so is finite. Consequently, $A(F^{\cyc})$
 is finitely generated. \epf

By the preceding lemma, it makes sense to speak of $\rank_{\Z}(A(F^{\cyc}))$ under the validity of conjecture of Mazur and Schneider. In view of Corollary \ref{rank formula corollary} and the question raised after, it seems plausible to make the following conjecture.

\bigskip\noindent\textbf{Conjecture 1.}
Assume that $(i)$ $A$ is an abelian variety over a number field $F$ which has ordinary reduction at every primes above $p$,
$(ii)$ $F_{\infty}$ is a uniform admissible $p$-adic extension of $F$ of dimension $d\geq 2$ and $(iii)$ $X(A/F^{\cyc})$ is finitely
generated torsion over $\Zp\ps{\Ga}$.
Then we have
\[ \rank_{\Z}(A(F_n)) \leq \left(\rank_{\Z}(A(F^{\cyc})) +
\sum_{\substack{w\in S(F^{\cyc}),
\\ \dim H_w\geq 1}}\corank_{\Zp}\big(Z_v(F^{\cyc}_w)(p)\big)\right)p^{(d-1)n} + O(p^{(d-2)n}),\]
where $S(F^{\cyc})$ is the set of primes of $F^{\cyc}$ above $S$ which are not divisible by $p$.
In the event that the extra assumption of Theorem \ref{main theorem2} is also valid, we replace $O(p^{(d-2)n})$ in the above inequality by $d\corank_{\Zp}(A(F_{\infty})(p))$.

\medskip
The point of Conjecture 1 is that here we have replaced $\rank_{\Zp}(X(A/F^{\cyc}))$ by the smaller quantity $\rank_{\Z}(A(F^{\cyc}))$. Note that there are examples where $\rank_{\Zp}(X(A/F^{\cyc})) \neq \rank_{\Z}(A(F^{\cyc}))$ (see \cite[pp. 140-142]{GR99}).
Nevertheless, we at least can record the following simple observation.

\bp \label{partial}
Assume that $(i)$ $A$ is an abelian variety over a number field $F$ which has ordinary reduction at every primes above $p$,
$(ii)$ $F_{\infty}$ is a uniform admissible $p$-adic extension of $F$ of dimension $d\geq 2$ and $(iii)$ $X(A/F^{\cyc})$ is finitely generated over $\Zp$ with $\rank_{\Zp}(X(A/F^{\cyc})) = \rank_{\Z}(A(F^{\cyc}))$.
Then we have
\[ \rank_{\Z}(A(F_n)) \leq \left(\rank_{\Z}(A(F^{\cyc})) +
\sum_{\substack{w\in S(F^{\cyc}),
\\ \dim H_w\geq 1}}\corank_{\Zp}\big(Z_v(F^{\cyc}_w)(p)\big)\right)p^{(d-1)n} + O(p^{(d-2)n}),\]
where $S(F^{\cyc})$ is the set of primes of $F^{\cyc}$ above $S$ which are not divisible by $p$.
\ep

\bpf
Since $X(A/F^{\cyc})$ is finitely generated over $\Zp$, it follows from a similar argument to that in \cite[Theorem 2.1]{CS12} that $X(A/F_{\infty})$ is finitely generated over $\Zp\ps{H}$ and hence satisfies the $\M_H(G)$-conjecture. The conclusion is now a consequence of Corollary \ref{rank formula corollary} and hypothesis (iii).
\epf

We also refer readers to \cite[Theorem 1.8]{DT},  \cite[Section 2.5]{DLMor} and \cite[Theorems A.38 and A.41]{DD} for discussion in support for Conjecture 1.

Finally, one may be tempted to ask the following naive question.

\medskip\noindent\textbf{Question.}
Retain the assumptions of Conjecture 1. Does one always have
\[ \rank_{\Z}(A(F_n)) = \left(\rank_{\Z}(A(F^{\cyc})) +
\sum_{\substack{w\in S(F^{\cyc}),
\\ \dim H_w\geq 1}}\corank_{\Zp}\big(Z_v(F^{\cyc}_w)(p)\big)\right)p^{(d-1)n} + O(p^{(d-2)n})\]
for $n\gg 0$?
\medskip

However, as we shall see in the next section, this question has a \textit{negative} answer.

\section{Additional evidence for Conjecture 1} \label{Additional evidence for Conjecture 2}

In this section, we describe how the deep works of Cornut-Vatsal \cite{CV}, Howard \cite{Howard} and Nekov\'{a}\u{r} \cite{Nek12} can be applied to give further evidence to Conjecture 1. For the remainder of this section, $E$ will denote an elliptic curve defined over $\Q$ with good ordinary reduction at the prime $p$. Let $F$ be an imaginary quadratic field of $\Q$, and $F_{\infty}$ the $\Zp^2$-extension of $F$. As before, write $G=\Gal(F_{\infty})$ and $H=\Gal(F_{\infty}/F^{\cyc})$. We also write $F_n$ for the intermediate subfield of $F_{\infty}/F$ with $\Gal(F_n/F)\cong\Z/p^n\times \Z/p^n$. We shall further assume that our elliptic curve $E$ has no complex multiplication and has conductor coprime to the discriminant of $F$. For ease of comparison, we first work out how the conjectured upper bound of Conjecture 1 looks like in this situation. As a start, we record the following lemma which we prove in slight generality.

\bl \label{finiteness}
Let $A$ be an abelian variety over $F$ with good ordinary reduction at the prime $v$ above $p$. Then for every prime $w$ of $F^{\cyc}$ above $v$, we have that $D_v(F^{\cyc}_w)$ is finite.
\el

\bpf
The long cohomology exact sequence of
\[ 0\lra C_v\lra A(p) \lra D_v\lra 0 \]
 gives rise to an exact sequence
 \[ A(F^{\cyc}_w)(p)\lra D_v(F^{\cyc}_w)\lra H^1(F^{\cyc}_w, C_v)\lra H^1(F^{\cyc}_w, A(p)) \lra H^1(F^{\cyc}_w, D_v) \lra 0,\]
 where the final zero follows from that fact that $H^2(F^{\cyc}_w,C_v)=0$ (cf. \cite[Theorem 7.1.8(i)]{NSW}). Since our abelian variety has good ordinary reduction at every prime of $v$, Imai's theorem \cite{Imai} asserts that $A(F^{\cyc}_w)(p)$ is finite. On the other hand, a local Euler characteristics argument (cf. \cite[\S 3]{Gr89}) shows that $\corank_{\Zp}(H^1(F^{\cyc}_w, A(p))) =
 \corank_{\Zp}(H^1(F^{\cyc}_w, C_v)) + \corank_{\Zp}(H^1(F^{\cyc}_w, D_v)$. Putting these information into the exact sequence, we see that the $D_v(F^{\cyc}_w)$ has trivial $\Zp$-corank.
\epf

By a theorem of Kato \cite{K}, $X(E/F^{\cyc})$ is torsion over $\Zp\ps{\Gal(F^{\cyc}/F)}$. Hence $\rank_{\Z}(E(F^{\cyc}))$ is well-defined (also see Theorem \ref{Kato Rohrlich}). Since $E$ is assumed to have no complex multiplication, $E(p)$ is not realizable over $F_{\infty}$ which, by \cite[Lemma 6.2]{LM} or \cite[Lemma 5.3]{Ze}, in turn implies that $E(F_{\infty})(p)$ is finite. Finally, since no primes outside $p$ ramified in a $\Zp^2$-extension, there are no extra contributions from local terms outside $p$. Hence Conjecture 1 in this setting reads as follow.

\bigskip\noindent\textbf{Conjecture $1'$.}
Retain the notation and settings of this section. Then
\[ \rank_{\Z}(E(F_n)) \leq \rank_{\Z}\big(E(F^{\cyc})\big)p^n.\]

\medskip
 Denote by $\epsilon(E/F, 1)$ the root number of Hasse-Weil $L$-function $L(E/F, s)$.

\bt[Cornut-Vatsal, Howard, Nekov\'{a}\u{r}] \label{Cornut-Vatsal}
Retain the above settings. In the event $\epsilon(E/F, 1)=-1$, suppose further that $F$ only contains one prime above $p$ and that $p$ does not divide the class number of $F$. Then
\[ \rank_{\Z}(E(F_n)) =\begin{cases}  O(1),& \mbox{if }~ \epsilon(E/F, 1)=+1 \\
      p^n+O(1), & \mbox{if }~ \epsilon(E/F, 1) =-1.\end{cases}\]
\et

\bpf
As explained in the proof of \cite[Proposition 3.14]{Van} (also see \cite{LV}), it follows from a combination of the deep results of Cornut-Vatsal \cite{CV}, Howard \cite{Howard} and Nekov\'{a}\u{r} \cite{Nek12} that
\[ \corank_{\Zp\ps{H}}(E(F_{\infty})\otimes \Qp/\Zp) =\begin{cases} 0, & \mbox{if }~ \epsilon(E/F, 1)=+1, \\
     1, & \mbox{if }~ \epsilon(E/F, 1) =-1.\end{cases}\]
Therefore, by an application of \cite[Theorem 1.10]{Har00}, we have
\[ \corank_{\Zp}\Big(\big(E(F_{\infty})\otimes \Qp/\Zp\big)^{H_n}\Big) =\begin{cases}  O(1),& \mbox{if }~ \epsilon(E/F, 1)=+1, \\
      p^n+O(1), & \mbox{if }~ \epsilon(E/F, 1) =-1.\end{cases} \]

On the other hand, we have
\[ \corank_{\Zp}(E(F_{n})\otimes \Qp/\Zp)  \leq \corank_{\Zp}\Big(\big(E(F_{\infty})\otimes \Qp/\Zp\big)^{G_n}\Big) \leq \corank_{\Zp}\Big(\big(E(F_{\infty})\otimes \Qp/\Zp\big)^{H_n}\Big),\]
where the second inequality is obvious (noting Lemma \ref{group}). For the first equality, we simply note that the kernel of the natural map
\[ E(F_{n})\otimes \Qp/\Zp\longrightarrow \big(E(F_{\infty})\otimes \Qp/\Zp\big)^{G_n} \]
is contained in the kernel of the map
\[ \Sel(E/F_n)\longrightarrow S(E/F_{\infty})^{G_n} \]
which in turn is contained in $H^1(G_n, E(F_{\infty})(p))$ as seen in the proof of Lemma \ref{main lemma}. But this latter group is finite as $E(F_{\infty})(p)$ is finite by the discussion before Conjecture $1'$.

Now the conclusion of the proposition clearly follows in the case of $\epsilon(E/F, 1)=+1$. For the situation when $\epsilon(E/F, 1)=-1$, Bertolini \cite[Proposition 7.6]{Be} has shown that
\[ \rank_{\Z}(E(F^{\mathrm{ac}}_n)) = p^n+O(1), \]
where here $F^{\mathrm{ac}}$ denotes the anticyclotomic $\Zp$-extension of $F$ and $F^{\mathrm{ac}}_n$ is the intermediate subfield of $F^{\mathrm{ac}}$ with $|F^{\mathrm{ac}}_n:F|=p^n$. Since $E(F^{\mathrm{ac}}_n)\subseteq E(F_n)$, the required equality follows.
    \epf

\br
Proposition 5.1, in particular, gives a negative answer to the question raised at the end of the preceding section when $\epsilon(E/F, 1)=+1$.
\er

\bc
 Retain the settings of Theorem \ref{Cornut-Vatsal}. In the event that $\epsilon(E/F, 1)=-1$, suppose further that $\Sha(E/F)(p)$ is finite. Then Conjecture $1'$ is valid.
\ec

\bpf
This is clear when $\epsilon(E/F, 1)=+1$. For the $\epsilon(E/F, 1)=-1$ case, it therefore remains to prove that $\rank_{\Z}(E(F^{\cyc}))\geq 1$ which amounts to showing that $\rank_{\Z}(E(F))\geq 1$. But by the parity result of \cite[Theorem 1]{Nek09}, and noting that $\Sha(E/F)(p)$ is finite by our hypothesis, we have that $\rank_{\Z}(E(F))$ is odd. In particularly, $\rank_{\Z}(E(F))\geq 1$.
\epf

\br
We emphasis that we \textit{do not} require the assumption that $X(E/F_{\infty})$ satisfies the $\M_H(G)$-conjecture throughout the discussion in this section.
\er

\section{A variant of Conjecture 1 for elliptic curve with supersingular reduction} \label{elliptic supersingular conjecture}

Throughout this section, let $E$ denote an elliptic curve over $\Q$ which has good supersingular reduction at the prime $p$. In particular, $E$ is no longer ordinary in the sense of Section 3. Despite this, we like to formulate a variant of Conjecture 1 for this class of elliptic curves in a modest setting, namely the case of a $\Zp^2$-extension.

Denote by $\widetilde{E}$ the reduced curve of $E$ modulo $p$. We shall assume that $a_p= p+1 - |\widetilde{E}(\mathbb{F}_p)|=0$ (note that this automatically holds if $p\geq 5$). Let $F$ be an imaginary quadratic field of $\Q$ at which the prime $p$ splits completely, say $p=\p\overline{\p}$. Denote by $F(\p^{\infty})$ the unique
$\Zp$-extension of $F$ unramified outside $\p$ and by $F(\p^n)$ the intermediate subfield of $F(\p^{\infty})$ with $|F(\p^n):F| =p^n$. We have analogous definitions for $F(\overline{\p}^{\infty})$ and $F(\overline{\p}^n)$. For each pair of nonnegative integers $m$ and $n$, write $F(\p^m\overline{\p}^n)$ for the compositum of the fields $F(\p^m)$ and $F(\overline{\p}^n)$.

We now denote by $\hat{E}$ the formal group associated to $E/\Qp$. Let $w$ be a prime of $F_{\infty}$ above $\p$. By abuse of notation, we write $w$ for the prime of $F(\p^m\overline{\p}^n)$ below this prime of $F_{\infty}$. Following \cite{Kim, LeiS}, we define the following groups
\[E^+(F(\p^m\overline{\p}^n)_w) = \{ P\in \hat{E}(F(\p^m\overline{\p}^n)_w)~:~\mathrm{tr}_{m/l+1,n}(P)\in \hat{E}(F(\p^l\overline{\p}^n)_w), 2\mid l, l<m \}, \]
\[E^-(F(\p^m\overline{\p}^n)_w) = \{ P\in \hat{E}(F(\p^m\overline{\p}^n)_w)~:~\mathrm{tr}_{m/l+1,n}(P)\in \hat{E}((F_{\p^l\overline{p}^n})_w), 2\nmid l, l<m \}, \]
where $\mathrm{tr}_{m/l+1,n}:\hat{E}(F(\p^m\overline{\p}^n)_w) \lra \hat{E}(F(\p^{l+1}\overline{\p}^n)_w)$ denotes the trace map.
For a prime $\overline{w}$ of $F_{\infty}$ above $\overline{\p}$, the groups
$E^{\pm}(F(\p^m\overline{\p}^n)_{\overline{w}})$ are defined in a similar fashion as above. From now on, write $F_n = F(\p^n\overline{\p}^n)$. Let $s,z\in\{+,-\}$. The signed Selmer group of $E$ over $F_n$ (cf. \cite{Kim}) is defined to be
\[\Sel^{s,z}(E/F_n)=\ker\left(\Sel(E/F_n)\lra \Bigg(\bigoplus_{w\mid \p}\frac{H^1(F_{n,w},E(p))}{E^s(F_{n,w})\ot\Qp/\Zp} \Bigg)\oplus \Bigg(\bigoplus_{\overline{w}\mid \overline{\p}}\frac{H^1(F_{n,\overline{w}},E(p))}{E^z(F_{n,\overline{w}})\ot\Qp/\Zp} \Bigg)\right).\]
Set $\Sel^{s,z}(E/F_{\infty})=\ilim_n \Sel^{s,z}(E/F_n)$, where $F_{\infty}$ is the $\Zp^2$-extension of $F$. We then write $X^{s,z}(E/F_{\infty})$ for the Pontryagin dual of $\Sel^{s,z}(E/F_{\infty})$.  We also write $G=\Gal(F_{\infty}/F)$ and $G_n=\Gal(F_{\infty}/F_n)$. Recently, Lei and Sprung have established the following result (see \cite[Theorem 4.4]{LeiS}).

\bt[Lei-Sprung]
Suppose that $X^{s,z}(E/F_{\infty})$ is torsion over $\Zp\ps{G}$ for every $s,z\in\{+,-\}$. Then one has $\rank_{\Z}(E(F_n))=O(p^n)$.
\et

It is then natural to ask if one can give an explicit upper bound as in the ordinary setting. This is the goal of the remainder of this section. Before doing so, we need to make the following supersingular analogue of the $\M_H(G)$-conjecture (also see \cite[Conjecture 3.16]{Lei} and \cite[Conjecture 5.3]{LeiZ}).

\medskip \noindent \textbf{Supersingular $\M_H(G)$-Conjecture.} \textit{For every $s,z\in\{+,-\}$, the module $X_f^{s,z}(E/F_{\infty})$ is finitely generated over $\Zp\ps{H}$, where $X_f^{s,z}(E/F_{\infty})= X^{s,z}(E/F_{\infty})/X_f^{s,z}(E/F_{\infty})(p)$.}

\medskip
The next result records an important consequence of the above conjecture on the structure of the module $X_f^{s,z}(E/F_{\infty})$ which we shall require later.

\bp \label{supersingular MHG} Retain the settings and notation of this section.
Assume further that the Supersingular $\M_H(G)$-Conjecture is valid for every $s,z\in\{+,-\}$.  Then  $H_i\big(H_n, X_f^{s,z}(E/F_{\infty})\big)=0$ for every $i\geq 1$ and $n\geq 0$.
\ep

\bpf
We first show that $X_f^{s,z}(E/F_{\infty})$ has no nonzero torsion $\Zp\ps{H}$-submodules.
By the Supersingular $\M_H(G)$-Conjecture, $X^{s,z}(E/F_{\infty})$ is in particular $\Zp\ps{G}$-torsion. Therefore, a similar argument to that in \cite[Theorem 1.1]{KimPN} can be applied to show that $X^{s,z}(E/F_{\infty})$ has no nonzero pseudo-null $\Zp\ps{G}$-submodules. By \cite[Lemma 4.2]{Su}, this in turn implies that $X_f^{s,z}(E/F_{\infty})$ has no nonzero pseudo-null $\Zp\ps{G}$-submodules.
Now a  well-known theorem of Venjakob \cite{V03} says that a $\Zp\ps{G}$-module $M$ which is $\Zp\ps{H}$-finitely generated is a pseudo-null $\Zp\ps{G}$-module if and only if it is a torsion $\Zp\ps{H}$-module. Since $X_f^{s,z}(E/F_{\infty})$ is finitely generated over $\Zp\ps{H}$ by the validity of the Supersingular $\M_H(G)$-Conjecture, the claim of this paragraph
then follows from a combination of these observations.

We now prove our proposition. Since $H_n\cong\Zp$, we have that $H_i(H_n,X_f^{s,z}(E/F_{\infty})) =  0$ for $i\geq 2$. Denoting by $\ga_H$ a topological generator of $H$, we have an identification $H_1(H_n,X_f^{s,z}(E/F_{\infty}))=X_f^{s,z}(E/F_{\infty})[\ga_H^{p^n}-1]$. But as seen in the previous paragraph, $X_f^{s,z}(E/F_{\infty})$ has no nonzero torsion $\Zp\ps{H}$-submodules. Therefore, we must have $H_1(H_n,X_f^{s,z}(E/F_{\infty})) =0$. The proof of the proposition is now complete.
\epf

We are in position to establish the following supersingular analogue of Theorem \ref{main theorem2}.

\bt \label{supersingular theorem}
Retain the settings and notation of this section. Assume further that Supersingular $\M_H(G)$-Conjecture is valid for all $s,z\in\{+,-\}$. Then we have
\[ \rank_{\Z}\big(E(F_n)\big)\leq \left(\sum_{s,z\in\{+,-\}}\rank_{\Zp\ps{H}}\big(X^{s,z}(E/F_{\infty})\big)\right)p^n.\]
\et

\bpf
By \cite[Proposition 4.3]{LeiS}, we have
\[ \corank_{\Zp}\big(\Sel(E/F_n)\big) \leq \sum_{s,z\in\{+,-\}}\corank_{\Zp}\big(\Sel^{s,z}(E/F_n)\big). \]
For each $s,z\in\{+,-\}$, the kernel of the natural map
\[ \Sel^{s,z}(E/F_n)\lra \Sel^{s,z}(E/F_{\infty})^{G_n}\]
is contained in the kernel of the map
\[ \Sel(E/F_n)\longrightarrow \Sel(E/F_{\infty})^{G_n} \]
which in turn is contained in $H^1(G_n, E(F_{\infty})(p))$ via a similar argument to that in the proof of Lemma \ref{main lemma}. We claim that this latter group is trivial.
Supposing for now the claim holds. Then we have the following inequality
\[ \corank_{\Zp}\big(\Sel^{s,z}(E/F_n)\big) \leq \rank_{\Zp}\big(X^{s,z}(E/F_{\infty})_{G_n}\big)\]
By Lemma \ref{rank compare lemma}, the term on the right is equal to $\rank_{\Zp}\big(X_f^{s,z}(E/F_{\infty})_{G_n}\big)$ which is less than or equal to $\rank_{\Zp}\big(X_f^{s,z}(E/F_{\infty})_{H_n}\big)$. By virtue of Proposition \ref{supersingular MHG}, we may apply Lemma \ref{Howson refined} to obtain $\rank_{\Zp}\big(X_f^{s,z}(E/F_{\infty})_{H_n}\big) = \rank_{\Zp\ps{H}}\big(X_f^{s,z}(E/F_{\infty})\big)p^n$. The conclusion of the theorem now follows from combining these observations.

It remains to verify our claim. To do this, it suffices to show that $E(F_{\infty})(p)=0$. Let $w$ be a prime of $F_{\infty}$ lying over $p$. Since $p$ splits completely over $F$, we have $F^{\cyc}_w = \Qp^{\cyc}$, and so we may apply \cite[Proposition 8.7]{Kob} to conclude that $E(F^{\cyc}_w)(p)=0$. It follows from this that $E(F^{\cyc})(p)=0$. Since $F_{\infty}/F^{\cyc}$ is a pro-$p$ extension, this in turn implies that $E(F_{\infty})(p)=0$ by \cite[Corollary 1.6.13]{NSW}.
\epf

We now proceed to formulate a variant of Conjecture 1 in this modest setting. Recall that it is well-known by now that $E$ is modular (see \cite{BCDT, Wi}). Therefore, one may apply the results of Kato \cite{K} and Rohrlich \cite{Ro} to conclude the following.

\bt[Kato, Rohrlich] \label{Kato Rohrlich}
Let $E$ be an elliptic curve over $\Q$ and $L$ a finite abelian extension of $\Q$. Then $E(L^{\cyc})$ is a finitely generated abelian group.
\et

In particular, the theorem implies that the quantity $\rank_{\Z}\big(E(F^{\cyc})\big)$ is well-defined.
Thus, we are in position to state our conjecture on the growth of Mordell-Weil ranks for an elliptic curve with supersingular reduction at $p$ in this modest setting.

\bigskip\noindent\textbf{Conjecture 2.}
 Let $E$ be an elliptic curve over $\Q$ which has good supersingular reduction at the prime $p$ with $a_p=0$. Denote by $F$ an imaginary quadratic field of $\Q$ at which the prime $p$ splits completely. Write $F_n$ for the intermediate subfield of the $\Zp^2$-extension $F_{\infty}$ of $F$ with $\Gal(F_n/F)\cong \Z/p^n\times\Z/p^n$. Then we have
\[ \rank_{\Z}(E(F_n))\leq 4\rank_{\Z}\big(E(F^{\cyc})\big)p^n.\]

\medskip
Although we are not able to relate $\rank_{\Zp\ps{H}}\big(X^{s,z}(E/F_{\infty})\big)$ to invariants coming from the cyclotomic level, we believe it should be related to certain cyclotomic invariants which bound the quantity $\rank_{\Z}\big(E(F^{\cyc})\big)$, henceforth the $``4"$ appearing in our conjecture. In view of Conjecture $1'$, one might even ask if the $``4"$ can be removed. We do not have an answer at this point of the writing.

\end{document}